\documentclass[11pt]{article}
\usepackage[margin=1in,vmargin=1in]{geometry}
\usepackage{tikz}
\usepackage{selectp}
\usepackage{graphicx}
\usepackage{url}
\usepackage{multirow}
\usepackage{setspace}
\usepackage{amssymb}
\usepackage{comment}
\usepackage{amsmath}
\usepackage{amsthm}

\newtheorem{conjecture}{Conjecture}

\title{A Conjecture about Periods in Subtraction Games}
\author{Mark Daniel Ward\footnote{Department of Statistics, Purdue University, 150 North
  University Street, West~Lafayette, IN~47907--2067 U.S.A.,
  Email: mdw@purdue.edu}}

\begin{document}
\maketitle
\begin{center}
\emph{Dedicated to Richard Guy, in celebration of his upcoming 100th birthday.}
\end{center}

\begin{abstract}
We make a conjecture that characterizes the
periods of the nim values in subtraction games with subtraction set of size 3.
\end{abstract}

\section{Introduction}

Richard~J. Nowakowski has, for many years, maintained a document of
unsolved problems in combinatorial games~\cite{GONC}; subtraction
games are the very first problem discussed.  Although they are
fundamentally important in the realm of combinatorial games, they have
been surprisingly challenging to completely analyze.

A subtraction game $S(s_1, s_2, s_3, \ldots)$ is played in much the
same way as Nim, but the number of beans that can be removed from a
heap during a player's turn is limited to the subtract set $\{s_1, s_2, s_3, \ldots\}$.

To follow the notation of Nowakowski, if the nim value of a heap of
size $h$ in a subtraction game is written as $n_h$, then the analogous
nim sequence for the subtraction game is $n_0n_1n_2n_3n_4\ldots$.
For many examples of these nim sequences for subtraction games,
see~\cite[pp.~84--85, Table~1]{WW2001}.
It is well known that, if the subtraction set has a finite size, then
the analogous nim sequence eventually is periodic, but the periods
of the nim values remain largely mysterious.

Mark Paulhus and Alex Fink have derived values of the
periods in two cases, for subtraction sets of size 3, namely,
in the case where $s_1 = 1$ and $s_2, s_3$ are arbitrary,
and in the case where $s_1 < s_2 < s_3 < 32$ (see~\cite{GONC}).
To the best of the author's knowledge, 
in the latter case $s_1 < s_2 < s_3 < 32$, Paulhus and Fink did not
have a formula that characterizes the periods, but rather,
they derived the values of the periods in these $\binom{32}{3}$
cases (without deriving a general formula to characterize the periods).

\section{Characterization of Periods}

In the present manuscript, we give a characterization of the possible
values of the periods in subtraction games with subtraction
sets of size~3.  We use ``gcd'' as shorthand for ``greatest common divisor.''

\begin{conjecture}
Consider a subtraction game $S(s_1,s_2,s_3)$ with $s_1 < s_2 < s_3$,
in which the nim sequence
$n_0n_1n_2n_3n_4\ldots$ eventually has period $p$.

\vskip 6pt
\noindent Case~I.  If $s_3 = s_1 + s_2$, define $0 \leq j < 2s_1$ so that $s_2 -
s_1 \equiv j \bmod 2s_1$.  Then the period $p$ can be precisely
characterized as follows:
$$
p = \begin{cases}
s_2+s_3-j& \hbox{if $0 \leq j < s_1$}, \\
(s_1)(s_2+s_3+j-2s_1)/\operatorname{gcd}(s_1,2s_1-j) & \hbox{if $s_1 \leq j < 2s_1$}.
\end{cases}
$$

\noindent Case~II.  If $s_3 \neq s_1 + s_2$, then $p$ has one of seven
potential values, which can be characterized as follows: 
The period $p$ is a divisor of \textbf{at least
one} of the numbers $s_i + s_j$ for $1 \leq i < j \leq 3$, and
moreover, the period $p$ is \textbf{exactly the gcd of all such terms}, i.e., 
$$p = \operatorname{gcd}\limits_{(i,j)\in\mathcal{G}}(s_i + s_j),$$
where $\mathcal{G}$ is the set of pairs $(i,j)$ such that $s_i+s_j$
is a multiple of $p$.
\end{conjecture}

The conjecture characterizes the period $p$ of every
subtraction game $S(s_1, s_2, s_3)$ with subtraction set of size~3.
I have verified the conjecture in all $\binom{4096}{3} =
11{,}444{,}858{,}880$ cases in which $1 \leq s_1 < s_2 < s_3 \leq 4096$.

\section{Acknowledgements}

\noindent\indent I first extend my gratitude and love to my wife Laura, in celebration
of 15 years of marriage.  

I am thankful for Doug Crabill's advice about the computations that
were used to help derive this conjecture.
I thank Kean Ming Tan and Xueyao Chen for early discussions about this
problem during an REU in summer 2008.

M.~D. Ward's research is supported by 
NSF Grant DMS-1246818, and by the NSF Science \& Technology Center for
Science of Information Grant CCF-0939370.

\bibliographystyle{plain}
\bibliography{subtraction}
\end{document}